\documentclass[12pt]{article}
\usepackage{amsmath}
\usepackage{amssymb}
\usepackage{amsthm}
\def\R{\mathbb{R}}

\def\N{\mathbb{N}}

\DeclareMathOperator{\exo}{exo}
\DeclareMathOperator{\card}{card}
\DeclareMathOperator{\Var}{{\mathbb Var}}
\def\qed{$\Box$}

\begin{document}

\numberwithin{equation}{section}

\newcommand{\I}{{\bf 1}}
\newcommand{\CF}{{\cal F}}
\newcommand{\calL}{{\cal L}}
\newcommand{\calN}{{\cal N}}

\newcommand{\calF}{{\cal F}}
\newcommand{\calB}{{\cal B}}
\newcommand{\calT}{{\cal T}}
\newcommand{\calE}{{\cal E}}
\newcommand{\calS}{{\cal S}}
\newcommand{\calA}{{\cal A}}
\newcommand{\calJ}{{\cal J}}
\newcommand{\calI}{{\cal I}}
\newcommand{\calK}{{\cal K}}
\newcommand{\uoclass}{{\calL}_{\rm uo}}
\newcommand{\loclass}{{\calL}_{\rm lo}}
\newcommand{\cclass}{{\calL}_{\rm c}}
\newcommand{\ismclass}{{\calL}_{\rm ism}}
\newcommand{\smclass}{{\calL}_{\rm sm}}
\newcommand{\ssmclass}{{\calL}_{\rm ssm}}
\newcommand{\cxclass}{{\calL}_{\rm cx}}
\newcommand{\idcxclass}{{\calL}_{\rm idcx}}
\newcommand{\dcxclass}{{\calL}_{\rm dcx}}
\newcommand{\stclass}{{\calL}_{\rm i}}
\newcommand{\lcxclass}{{\calL}_{\rm lcx}}
\newcommand{\plcxclass}{{\calL}_{\rm plcx}}
\newcommand{\icxclass}{{\calL}_{\rm icx}}
\newcommand{\icvclass}{{\calL}_{\rm icv}}
\newcommand{\cvclass}{{\calL}_{\rm cv}}
\newcommand{\extclass}{{\cal E}_{{\rm c}}}
\newcommand{\extclassc}[1]{{\cal E}_{{\rm {#1}}}}
\newcommand{\plcxcxclass}{{\calL}_{{\rm plcx}}}
\newcommand{\veco}[3]{\left({#1}_{#2},\ldots,{#1}_{#3}\right)}
\newcommand{\supp}{{\rm supp}}
\newcommand{\eqdistr}{\stackrel{\rm d}{=}}
\newcommand{\convdist}{\stackrel{\rm d}{\Longrightarrow}}
\newcommand{\convdistr}{\stackrel{\rm d}{\Longrightarrow}}
\newcommand{\convas}{\stackrel{\rm a.s.}{\Longrightarrow}}
\newcommand{\convprob}{\stackrel{p}{\Longrightarrow}}
\newcommand{\conv}[1]{\stackrel{{#1}}{\Longrightarrow}}
\newcommand{\ordersm}{<_{\rm sm}}
\newcommand{\orderism}{<_{\rm ism}}
\newcommand{\orderssm}{<_{\rm ssm}}
\newcommand{\ordercx}{<_{\rm cx}}
\newcommand{\ordericx}{<_{\rm icx}}
\newcommand{\orderlcx}{<_{\rm lcx}}
\newcommand{\orderplcx}{<_{{\rm plcx}}}
\newcommand{\orderst}{<_{\rm st}}
\newcommand{\orderidcx}{<_{\rm idcx}}
\newcommand{\ordercv}{<_{\rm cv}}
\newcommand{\ordersl}{<_{\rm sl}}
\newcommand{\orderse}{<_{\rm se}}
\newcommand{\orderdcx}{<_{\rm dcx}}
\newcommand{\orderconc}{<_{\rm c}}
\newcommand{\orderuo}{<_{\rm uo}}
\newcommand{\orderlo}{<_{\rm lo}}
\newcommand{\orderodwlo}{>_{\rm lo}}
\newcommand{\orderccx}{<_{\rm ccx}}
\newcommand{\ordericcx}{<_{\rm iccx}}
\newcommand{\ordersymcx}{<_{\rm symcx}}

\renewcommand{\theequation}{\mbox{\arabic{section}.\arabic{equation}}}
\newtheorem{proposition}{Proposition}[section]
\newtheorem{theorem}[proposition]{Theorem}
\newtheorem{corollary}[proposition]{Corollary}
\newtheorem{lemma}[proposition]{Lemma}
\newtheorem{definition}[proposition]{Definition}
\newtheorem{remark}[proposition]{Remark}
\newtheorem{example}[proposition]{Example}
\newtheorem{assumption}[proposition]{Assumption}
\newtheorem{fig}[proposition]{Figure}

\newcommand{\comments}[1]{\marginpar{ \begin{minipage}{0.5in} \tiny #1 \end{minipage} }}

\title{Comparisons and asymptotics for empty space
hazard functions of germ-grain models}
\author{G\"unter Last
\thanks{address: Institut f\"ur Stochastik, 
KIT-Campus S\"ud D-76128, Germany.}
\\
{\it  Karlsruher Institut f\"ur Technologie}
\and Ryszard Szekli
\thanks{ Work supported by the
    Marie Curie Fellowship for the Transfer of
    Knowledge Harmonic Analysis, Nonlinear Analysis and Probability,
    MTKD-CT-2004-013389, and MENiS-KBN 2007-2009 grant  N N2014079 33, address: 
Mathematical
    Institute,  pl. Grunwaldzki 2/4,  50-384 Wroc{\l}aw, Poland.}
 \\
  {\it Wroc{\l}aw University}
} \maketitle

\begin{abstract} We study stochastic properties
of the empty space for stationary germ-grain models in $\R^d$, in
particular we deal with the inner radius of the empty space with
respect to a general structuring element which is allowed to be
lower-dimensional.  We consider Poisson cluster germ-grain models
and Boolean models with grains that are clusters of convex bodies
and show that more variable size of clusters results in
stochastically greater empty space in terms of the empty space
hazard function. We also study impact of clusters being more
spread in the space on the value of the empty space hazard. Further we
obtain asymptotic behavior of the empty space hazard
functions at zero and at infinity.
\end{abstract}



\section{Introduction}
The statistical analysis of a spatial pattern $Z\subset\R^d$ is based on the
assumption that $Z$ is a random set in $\R^d$. Distance methods for point patterns usually begin by estimating nonparametric summary functions such as the empty space function, nearest neighbor distance distribution, Ripley's $K$-function and derived 
statistics such as the pair correlation and the $J$-function. 
For a spatial pattern $Z\subset\R^d$, which is not necessarily a point process, 
such as germ-grain model a particularly
useful functional for estimating properties of $Z$ is the empty 
space hazard (rate) function. Smaller empty space hazard function in a 
germ-grain pattern with the same germ intensity intuitively corresponds 
to a more clustered pattern with more empty space. In this paper we shall 
give some sufficient conditions for the hazard rate ordering of the empty 
space distribution  in two compared germ-grain models using particular 
classes of such spatial patterns such as Poisson cluster point processes, 
Poisson cluster germ-grain models, and mixed Poisson germ-grain models. 
Our sufficient conditions will be given in a language of stochastic 
orderings which we define for some parameters of these models. Some early 
results on stochastic ordering of random closed sets can be found in
Stoyan and Stoyan \cite{SS80}, and a characterization of the strong stochastic 
ordering of random closed sets is given by Molchanov \cite{Mo05} in Theorem 4.42.
However, apart from Section 3.8 of Hall \cite{Hall88} (dealing
with volume fractions) and Last and Holtmann \cite{LH99}
(dealing with the spherical contact distribution of a Gauss-Poisson model) 
we are not aware of papers comparing functionals of stationary random closed sets.
We introduce two stochastic ordering relations which  are  useful for comparison 
two distributions with equal expected values in which case both 
relations imply that a larger variable in these orderings has a 
bigger variance. It turns out  that the impact of a larger variance 
in a spatial pattern for some  component of a germ grain model with 
fixed expectation is that a larger variance  gives more empty space with 
more clustering.  But the situation is more complicated if we compare 
models with equal intensities of Poisson cluster germ points but at the same time 
increasing or decreasing intensities of underlying Poisson processes.

To be more precise, for two random variables $\eta$, $\tilde\eta$, 
taking values $0,1,2,\ldots$ we say that they 
are ordered in the length biased probability generating functions ordering, and 
write $\eta<_{l-g}\tilde\eta$ if the corresponding length biased 
variables $\eta_l$ and $\tilde\eta_l$ are ordered in the 
probability generating ordering, see \cite[Section 1.8]{St83}.  
For two nonnegative random variables $\Lambda_1$, $\Lambda_2$ 
we say that they are ordered in the first cumulant order and 
write $\Lambda<_{cum}\tilde\Lambda$ if for the corresponding 
cumulant generating functions $C_{\Lambda}$,  $C_{\tilde\Lambda}$, 
taking the first derivatives we have  $C'_{\Lambda}(s)\ge C'_{\tilde\Lambda}(s)$, 
$s\in [-1,0]$. We show for  Neyman-Scott processes with cluster sizes 
$\eta$, $\tilde\eta$ that
$\eta<_{l-g}\tilde\eta$ implies that the corresponding empty space distributions 
are ordered in the hazard rate ordering. For mixed Poisson germ-grain models 
we show that if the random intensities are ordered $\Lambda<_{cum}\tilde\Lambda$ 
then the corresponding empty space distributions are ordered in the hazard 
rate ordering. For both  introduced orderings, under the additional assumption 
that the ordered variables have equal expected values, the greater 
variable has larger variance. Therefore the mixed Poisson germ-grain 
model behaves similarly as the Neyman-Scott germ-grain model with respect to 
variance changes inside the model. We study also behavior of the 
Gauss-Poisson model. For details, additional properties of the orderings, and 
examples see Sections 4 and 5. We shall also study there the asymptotic behavior of
empty space hazard rates  at zero and at infinity.

\section{Preliminaries}

\subsection{Empty space hazard functions}
Let us recall the definition of the empty space distribution function $F$ (called also the
first contact distribution function). For $\|x\|$ the
Euclidean norm of a vector $x\in\R^d$, and $
d(x,A):=\inf\{\|x-y\|:y\in A\} $ the distance between $x\in\R^d$
and a set $A\subset \R^d$, this distribution is given by
$$
F(t)=P(d(0,Z)\le t),\quad t\ge 0,
$$
where $0$ denotes the zero vector. The value $F(0)=P(0\in Z)$ is
the {\em volume fraction} of $Z$. $F$ can also be written in terms
of the {\em capacity functional} $T_Z$ (defined by
$T_Z(K):=P(Z\cap K\ne \emptyset),\quad K\subset\R^d\;{\rm
compact}$), namely $F(t)=P(Z\cap B(x,t)\ne
\emptyset)=T_Z(B(x,t))$, $t\ge 0$, where $B(x,t)$ is the closed
ball with center $x$ and radius $t$. Stationarity of $Z$ ensures
that $F(t)$ does not depend on $x$. Hence, $F(t)$ is the
probability that $Z$ hits the ball $B(t):=B(0,t)$. There are many
reasons for studying other distances than the Euclidean distance.
For example digital image analysers estimate rather polygonal
distance than spherical one. To estimate isotropy of point
patterns one needs elliptical distance. It is clear that the
distribution of $Z$ is not determined by $T_Z(B(x,t))$ for all
balls, and a larger class of sets than the class of balls provides
a better information on $Z$.

The usual way of introducing other distances than the spherical one is
to fix a {\em structuring element} ({\em gauge body})
$B\subset\R^d$. This is a compact convex set having $0\in B$. Then
the $B$-distance of a point $x\in{\R^d}$ to a set $K\subset\R^d$
is defined by
$$
d_B(x,K):=\inf\{t\ge 0: (x+tB)\cap K\ne \emptyset\}.
$$
It is possible that the set on the right side is empty, e.g.\ if
$B$ is lower-dimensional. In such a case we set
$d_B(x,K):=\infty$. Note that we have the translation invariance
property  $d_B(z+x,z+K) = d_B(x,K)$, for all $z\in\R^d$. Clearly,
$d_B(x,K)\le t$ if and only if $x$ is contained in the {\em
generalized outer parallel set} $K+ tB^\ast$ of $K$, where $B^*$
denotes the reflected set $\{-x:x\in B\}$. If $B$ is full
dimensional (i.e.\ has a non-empty interior) and centrally
symmetric (i.e.\ $B^*=B$), then $d_B(\cdot,\cdot)$ is a metric on
$\R^d$ induced by the norm $d_B(\cdot,0)$, and the pair
$(\R^d,d_B)$ is called a {\em Minkowski space}.

If the $B$-distance $d_B(x,K)$ of a point $x\notin K$ is attained
in a unique point $y$ in the boundary $\partial K$ of $K$ (that
means, if $(x+d_B(x,K) B)\cap K =\{ y\}$), then we define the {\em
relative metric projection} of $x$ on $K$ by $p_B(K,x):=y$, and
the {\em contact direction vector} $u_B(K,x)$ as the element of
$\partial B^*$ given by
$$
u_B(K,x):=\frac{x-y}{d_B(x,K)}.
$$
The points  $x\in\R^d\setminus K$ for which the distance $d_B(x,K)$ is attained
in more than one point of $K$ (and for which $u_B(x,K)$ is therefore not
defined) form the {\em exoskeleton} $\exo_B(K)$ of $K$ (see Hug, Last and Weil
\cite{HLW02a}). In the Euclidean case, and if $K$ is a finite or locally finite
set, $\exo_{B(1)}(K)$ is (the  boundary of) the {\em Voronoi tessellation}
generated by $K$.

We define the {\em directed, $B$-relative empty space function}
$F_B$ of $Z$ by
\begin{align}\label{cdf2}
F_B(t,C):=P(d_B(0,Z)\le t,u_B(Z,0)\in C),\quad t\ge 0,\,C\in\mathcal{B}^d,
\end{align}
where $\mathcal{B}^d$ is the system of Borel subsets of $\R^d$.
Here we use the
convention $u_B(Z,0):=u_0$ if $0\in Z$ or $d_B(0,Z)=\infty$, where
$u_0\in\partial B^*$ is fixed. Definition \eqref{cdf2} is subject
to the assumption that the vector $u_B(Z,0)$ is $P$-a.s.\
well-defined on $\{0<d(0,Z)<\infty\}$. If $Z$ is a random closed
set and $B$ is strictly convex, containing $0$ in its interior,
then this is indeed the case. This follows from the fact
that $\exo_B(Z)$ has volume $0$ and from stationarity of $Z$.
(More general cases require a
suitable assumption on the relative positions of $Z$ and $B$,
see Subsection \ref{subgermgrain}).
The  function $F_B$ determines the
joint distribution of the pair $(d_B(0,Z),u_B(Z,0))$, and hence
that of the {\em contact vector} $d_B(0,Z)u_B(Z,0)$. For each
fixed $t$, the function $F_B(t,\cdot)$ is a measure on $\R^d$
concentrated on $\partial B^*$. The function
$F_B(\cdot):=F_B(\cdot,\R^d)$ is called ({\em $B$-relative) empty
space function} of $Z$.

The {\em Minkowski addition} of two sets $C,D\subset \R^d$
($C\oplus D:=\{x+y:x\in C,y\in D\}$) gives another form for $F_B$.
Stationarity easily implies that
\begin{align}\label{2.3a}
F_B(t):=F_B(t,\R^d)= V_d(A)^{-1} E[ V_d((Z\oplus tB^*)\cap A)],\quad t\ge 0,
\end{align}
for each Borel test set $A$, such that the {\em volume}
$V_d(A)$ of $A$ is positive and finite, see e.g.\ \cite{HLW02b,BG99}.

Hansen, Baddeley and Gill \cite{BG99}, utilizing Federer's coarea
theorem, showed that the empty space function $F_B$ of a random
closed set $Z$ is absolutely continuous on $(0,\infty)$ with
density
\begin{align}\label{bdensity}
f_B(t)= V_d(A)^{-1}E\left[\int _{A\cap \partial(Z\oplus tB^*)}\|\nabla
d_B(s,Z)\|^{-1}\mathcal{H}^{d-1}(ds)\right]
\end{align}
where $\mathcal{H}^i$, $i\in\{0,\ldots,d\}$, denotes
$i$-dimensional Hausdorff measure on $\R^d$, $\partial A$ denotes the
boundary of  $A$, and  $\nabla d_B$ denotes the gradient of the
function $d_B$. In the Euclidean case, this formula reduces to
\begin{align}\label{density}
f(t)=V_d(A)^{-1} E[\mathcal{H}^{d-1}(\partial(Z\oplus B(t))\cap A)],
\end{align}
In this case the empty space hazard equals the ratio of the
expected measure of the boundary $\partial(Z\oplus B(t))$ inside
the window $A$ to the volume of of the space not occupied by
$Z\oplus B(t)$ inside $A$.

For general $B$ the empty space hazard is given  by
\begin{align}\label{brate}
r_B(t)= \frac{1}{V_d(A)-E[ V_d(Z\oplus tB^*)\cap A)]}
E\left[\int_{A\cap
\partial(Z\oplus tB^*)}\|\nabla d_B(s,Z)\|^{-1} \mathcal{H}^{d-1}(ds)\right],
\end{align}
which is intuitive in the sense that the empty space hazard
depends on the speed of increase of the distance function $d_B$
along all coordinates.
%\textbf{can we see from this formula how the shape of $B$
%influences the empty space**}
It is possible to rewrite this formula in terms of the support
function $h_B$ of the gauge body $B$ (see \cite{BG99}).
\begin{align}\label{bratesupp}
r_B(t)= \frac{1}{V_d(A)-E[V_d(Z\oplus tB^*)\cap A)]}
E\left[\int_{A\cap\partial(Z\oplus tB^*)}(h_B(u_B(Z,s)))^{-1}
\mathcal{H}^{d-1}(ds)\right].
\end{align}

The direction dependent (sub)distribution functions,  $F_B(\cdot,C)$, as defined in
(\ref{cdf2}) are also absolutely continuous on $(0,\infty)$ for
any $C\in \mathcal{B}^d$.
%\textbf{can we cite something for abs. cont. of} $F_B(\cdot,C)$?
Letting $f_B(\cdot,C)$ denote its density, we define
\begin{align}\label{emptyhazard}
r_B(t,C):=\frac{f_B(t,C)}{1-F_B(t)},
\end{align}
where $a/0:=0$ for all $a\in\R$. We call the function
$r_B(\cdot,C)$ the directed, {\em $B$-relative empty space hazard
of $Z$}.

 For fixed  $B$,  we shall order two random sets $Z$, and $\tilde Z$ with 
respect to their $B$-relative empty space hazard functions, and write
\begin{align}\label{2.8}
Z<_{h-B} \tilde Z
\end{align}
iff for all $t\ge 0$, and $C\in\mathcal{B}^d$
\begin{align}\label{monhazard8}
r_B(t,C)\ge \tilde r_B(t,C),\quad t\ge 0.
\end{align}
This ordering is stronger than the usual (strong) stochastic ordering of empty 
space distributions, for further details on such orderings see e.g.\ 
Szekli \cite[Section 1.4]{Szekli95}.

\subsection{Empty space hazard rates via support measures}\label{subgermgrain}

By a  {\em germ-grain model} in $\R^d$ we mean a random set of the form
$$
Z=\bigcup^\infty_{n=1}(X_n +\xi_n) =\bigcup^\infty_{n=1}\{x
+\xi_n:x\in X_n\},
$$
where the random points $\xi_n$, $n\in\N$, represent the locations
of the {\em germs} and the {\em primary grains} $X_n$, $n\in\N$,
are assumed to be random non-empty compact subsets of $\R^d$. We
assume that the (simple) {\em point process} $N:=\{\xi_n:
n\in\N\}$ is {\em stationary}, that is the distribution of the
{\em shifted} point process $N+x:=\{\xi_n+x:n\in\N\}$ does not
depend on $x\in\R^d$, and  that $N$ is independent of $(X_n)_{n\ge
1}$ which is a sequence of independent and identically
distributed random sets. The {\em intensity}
$\lambda:=E\card\{n\in\N:\xi_n\in[0,1]^d\}$ of $N$ is assumed to
be finite. An important special case is the {\em Boolean model},
where the germs are located according to a homogenous {\em Poisson
process}. The grains $X_n$ as well as the germ-grain model $Z$
itself are measurable mappings from $\Omega$ into the set
$\mathbf{F}$ of all closed subsets of $\R^d$. Measurability refers  to the
smallest $\sigma$-field of subsets of $\mathbf{F}$, containing the
sets $\mathbf{F}_K=\{F:F\cap K\ne \emptyset\}$, for all compact sets
$K\subset \R^d$. Stationarity of $N$ entails that also $Z$ is
stationary, i.e.\ that the distribution of $Z+x$ does not depend
on $x$.
It is convenient to denote by $X_0$  a {\em typical grain}
having its distribution equal to that of $X_n$.
We assume that $E[V_d(X_0+K)]$ is finite
for all compact $K\subset\R^d$.
We will use later that the capacity functional of a
Boolean model is given by
\begin{align}\label{CFBoolean}
P(Z\cap K\ne\emptyset)=1-\exp[-\lambda E[V_d(X_0+K^*)]].
\end{align}
In particular, the  volume fraction of a Boolean model is given by
\begin{align}\label{vfBoolean}
P(0\in Z)=1-\exp[-\lambda E[V_d(X_0)]].
\end{align}
We refer to Stoyan, Kendall and Mecke \cite{SKM95} and Schneider and Weil
\cite{SW08} for a detailed introduction to germ-grain models.

Consider a convex, compact and non-empty set $K\subset \R^d$. We
assume that $K$ and $B^*$ are in {\em general relative position},
which means that $K$ and $B^*$ do not contain parallel segments in
parallel and equally oriented support (hyper)planes.
This means that $K$ and $B$ have {\em independent support}
sets, see \cite[p.\ 611]{SW08} for more detail.
A sufficient condition is that $K$ or
$B$ is strictly convex, This assumption guarantees that $p_B(K,x)$
(and hence $u_B(K,x)$) is defined for all $x\notin K$. Then there
are finite measures $C_0(K;B;\cdot),\ldots,C_{d-1}(K;B;\cdot)$ on
$\R^d\times\R^d$ which satisfy the {\em local Steiner formula}
\begin{multline}\label{2.1}
V_d(\{x\in\R^d: 0<d_B(x,K)\le t, (p_B(K,x),u_B(K,x))\in A\times C\})\\
=\sum^{d-1}_{i=0}t^{d-i}b_{d-i}C_i(K;B;A\times C)
\end{multline}
for all  $A,C\in\mathcal{B}^d$, where $b_i$
($i\in\N$) denotes  the volume of the unit ball in $\R^i$, and
$b_0:=1$. These {\em relative support measures} of $K$  are
uniquely determined by \eqref{2.1}. They are concentrated on
$\partial K\times\partial B^*$ and in fact on the {\em relative
normal bundle}
$$
N_B(K):=\{(p_B(K,x),u_B(K,x)): x\notin K\}
$$
of $K$. If $B=B(1)$ then the measures
$C_i(K;\cdot):=C_i(K;B(1);\cdot)$ are the {\em generalized
curvature measures} of $K$. The total mass
$V_i(K):=C_i(K;\R^d\times\R^d)$ is the $i$th {\it intrinsic
volume} of $K$. In particular, $V_d(K)$ is the volume of $K$,
$V_{d-1}(K)$ equals one half of the surface area, $V_{d-2}(K)$ is
proportional to the integral mean curvature, $V_1(K)$ is
proportional to the mean width of $K$, and $V_0(K)=1$. Equation
\eqref{2.1} implies the classical {\em Steiner formula}
\begin{align} \label{Steiner}
V_d(K\oplus B(t))=\sum^d_{i=0}b_{d-i}t^{d-i}V_i(K).
\end{align}
In the general case, the total mass $C_i(K;B;\R^d\times\R^d)$
is a special {\em mixed volume}, namely
\begin{align}\label{mixed}
C_i(K;B;\R^d\times\R^d)=b_{d-i}^{-1}\binom{d}{i}V(K[i],B^*[d-i]),
\quad i=0,\ldots,d-1.
\end{align}
For $i=0$ we have
\begin{align}\label{mixedzero}
C_0(K;B;\R^d\times\R^d)=b_{d}^{-1}V_d(B^*),
\end{align}
% (\textbf{$C_0$ here is not dependent on $K$, is it also true for
% $C_0(K;B;A\times C)$** })
%and for $B=B(1)$
%\begin{multline}\label{mixedzeroball}
%C_0(K;B(1);A\times C)\\
%=(db_{d})^{-1}\int_{S^{d-1}}\I\{u\in C,
%u_{B(1)}(K,x)=u, p_{B(1)}(K,x)\in A, x\notin K\}\mathcal{H}^{d-1}(du),
%\end{multline}
see \cite{SW08} for the notation used here and
for further details on support and curvature measures.
%\textbf{More details for Schneider and Weil, pages** }

Consider now  a germ-grain model $Z$ with convex, compact
grains.  Assume that the {\em reduced second moment measure} of
$N$ on $\R^d$, defined by
\begin{align}\label{reduced}
E\Bigg[\sum_{\substack{x,y\in N\\x\ne y}}
\I\{x\in[0,1]^d\}\{x-y\in\cdot\}\Bigg].
\end{align}
is absolutely continuous, and assume that the typical grain $X_0$
and $B^*$ are a.s.\ in general relative position. It then follows,
that $u_B(Z,0)$ is almost
surely well-defined on $\{0<d_B(0,Z)<\infty\}$. This can be
proved as Proposition 4.9 in \cite{Hl00}.

By the Steiner formula \eqref{Steiner}
our general integrability assumption on $X_0$ (see
Subsection \ref{subgermgrain}) is equivalent to
the finiteness of the mean intrinsic volumes
\begin{equation}\label{moment}
\bar V_i:=E V_i(X_0),\quad i=0,\ldots,d,
\end{equation}
of the typical grain. The Steiner formula \eqref{Steiner}
together with the local Steiner formula \eqref{2.1} imply that
\begin{equation}\label{relmoment}
\bar V_{i,B}:=E C_i(X_0;B;\R^d\times\R^d),\quad i=0,\ldots,d-1,
\end{equation}
are finite as well. Therefore the mean relative support measures of the
typical grain, defined by
$$
\bar C_{i,B}(\cdot):=E C_i(X_0;B;\cdot),\quad i=0,\ldots,d-1,
$$
are finite measures on $\R^d\times \R^d$.

We further use the {\em Palm probability}  $P^0_N$ of $P$ with
respect to $N$ (see \cite{DaVJ03,SKM95}). We
can interpret $P^0_N(A)$ as the conditional probability of the
event $A\in\CF$ given that $0$ is a ``randomly chosen point'' of
$N$.  Let us define $X(x):=X_n$ if $x=\xi_n$ for some $n$, and
$X(x):=\emptyset$, otherwise.  Then under  $P^0_N$,
$\{(x,X(x)):x\in N\}$ is an independently marked point process,
i.e.\ conditionally on $N$, the  grains $\{X(x):x\in N\}$
are independent and have the same distribution as $X_0$.
For $\bar V_{i,B}>0$, $t\ge 0$, and $C\in \mathcal{B}^d$, let
\begin{align}\label{GiB}
G_{i,B}(t,C):= \bar V_{i,B}^{-1}\int \I\{u\in C\}P^0_N(d_B(y+tu,Z^{!})\le t)
\bar C_{i,B}(d(y,u)),
\end{align}
where
\begin{align}\label{Z!}
Z^{!}:=\bigcup_{x\in N\setminus \{0\}}(X(x)+x)
\end{align}
is the union of all grains except for the  grain located at the
origin. We set $G_{i,B}\equiv 0$ for $\bar V_{i,B}=0$. The
function $G_{i,B}(\cdot,\R^d)$ can be interpreted as the
distribution function of a random variable $\xi$, say, which can
be constructed as follows. First one selects a  point $Y$ of $N$
at random. Then one samples a random element $(X,W)$  according to
the distribution $\bar{V}^{-1}_{i,B}\bar C_{i,B}$.
%\textbf{(I do not see this construction, which $i$ do we select**)}
If $Y$ is not covered by $\bigcup_{x\in N\setminus\{Y\}}(X(x)+x)$,
then $\xi$ is the $B$-distance from $Y+X$ to the exoskeleton
$\exo_B(Z)$, in the direction $W$. Otherwise $\xi=0$. We shall utilize
the following functions
\begin{align}\label{JiB}
J_{i,B}(t,C):=\frac{G_{i,B}(\infty,C)-G_{i,B}(t,C)}{1-F_B(t)},
\quad i=0,\ldots,d-1,
\end{align}
where $i\in\{0,\ldots,d-1\}$, $t\ge 0$, $C\subset\R^d\times\R^d$
is a Borel set, and
$$
G_{i,B}(\infty,C):=\lim_{t\to\infty}G_{i,B}(t,C)=
\bar V_{i,B}^{-1}\bar C_{i,B}(\R^d\times C).
$$
Special cases of these functions were introduced in
\cite[Section 5]{LH99} after the point process case
had been treated in \cite{LiBa96}.
The functions $J_{i,B}(t,C)$ can be used as non-parametric
measures for expressing differences between a general germ-grain
model and Boolean model with the same values of $\lambda
\bar{V}_{i,B}$. Intuitively speaking, such measures detect
interactions and clustering effects. In the
Euclidean case (and for $C=\R^d$) the following theorem was proved
in \cite{LaScha89,LH99}.
The following relative version
is implicit in Hug and Last \cite{Hl00}, at least in the case of a
strictly convex $B$. The general result can be derived from
Theorem 5.1  in Hug, Last and Weil \cite{HLW02a}. 

\begin{theorem} \label{lastholtmann} Consider a stationary
germ-grain model satisfying the assumptions formulated above.
Then $F_B$ is absolutely continuous and the $B$-relative
empty space hazard $r_B$ is given by
\begin{equation}\label{dens1}
r_B(t,C)=\sum^{d-1}_{i=0}(d-i)t^{d-i-1}b_{d-i}\lambda\bar{V}_{i,B}
J_{i,B}(t,C).
\end{equation}
\end{theorem}

If $N$ is a Poisson process (i.e.\ $Z$ is the Boolean model with convex grains),
then Slivnyak's theorem (see e.g.\ Stoyan, Kendall and Mecke \cite{SKM95})
implies that
\begin{equation}
P^0_N (d_B(y+tu,Z^{!})>t)=1-F_B(t).
\end{equation}
Hence \eqref{dens1} simplifies to
\begin{equation}\label{densBoolean}
r_B(t,C)=\sum^{d-1}_{i=0}(d-i)b_{d-i}t^{d-i-1}\lambda\bar S_{i,B}(C),
\end{equation}
where $\bar S_{i,B}(C):=\bar C_{i,B}(\R^d\times C)$. In the case
of a strictly convex gauge body $B$ this result can be found in
\cite{Hl00}. Note that in the Boolean model
$r_B(\cdot,\R^d)$ is determined by the intensity $\lambda$ and the
mean mixed volumes $\bar{V}_{1,B},\ldots, \bar{V}_{d-1,B}$ of
$X_0$. For $d=2$ and $B=B(1)$, for instance, the only parameter of
$X_0$ influencing the empty space hazard rate is its mean boundary
length $\bar V_1$. Note that in the Boolean model with convex
grains the empty space hazard rate is increasing, and
asymptotically
$$
\lim_{t\to \infty}t^{1-d}r_B(t,C)=\lambda d b_d
\bar S_{0,B}(C)=\lambda d b_d E [C_0(X_0;B;\R^d\times C)].
$$
By \eqref{mixedzero} we obtain in particular that
$$
\lim_{t\to \infty}t^{1-d}r_B(t,\R^d)=\lambda d b_d
E[b_{d}^{-1}V_d(B^*) ]=\lambda d V_d(B^*).
$$

%\textbf{Does $E [C_0(X_0;B;\R^d\times C)]$ depend on $X_0$** }

%The following result can be proved as in Last and Holtmann \cite{LH99}.
%
%\begin{theorem} Assume that $Z=Z^{(1)}\cup Z^{(2)}$
%is the superposition of two germ-grain models
%\[
%Z^{(k)}=\bigcup^\infty_{n=1}(X_{n,k} +\xi_{n,k}),
%\]
%satisfying our assumptions and such that the marked point processes
%$\{(\xi_{n,k},X_{n,k}):n\ge 1\}$, $k=1,2$, are independent. For $k=1,2$, let
%$\lambda_k$ denote the intensity of $N_k:=\{\xi_{n,k}:n\ge 1\}$ and
%$J_{i,B,k}$ the $J_{i,B}$-function of $Z^{(k)}$. Then, for $i=0,\ldots,d-1$,
%$$
%J_{i,B}(t)=\frac{\lambda_1}{\lambda_1+\lambda_2}J_{i,B,1}(t)
%+\frac{\lambda_2}{\lambda_1+\lambda_2}J_{i,B,2}(t),
%$$
%provided that $\bar V_i>0$.
%\end{theorem}

\section{Results for Poisson cluster point processes}

In this section we shall consider germ-grain models where the germ
process will be a Poisson cluster point process, and grains will
be one point grains attached to germs, hence $Z=N$. Such a model
is a special case of a germ-grain model with convex (one point)
grains but, alternatively, this model might be seen as a Boolean
model (germs form a Poisson process) with non-convex grains (point
clusters). Recall from \cite{DaVJ03} that a Poisson  cluster point
process can be written as
\begin{equation}
N=\bigcup_{x\in \Pi}L_x +x,
\end{equation}
where $\Pi$ is a Poisson process with positive and finite intensity
$\lambda_\Pi$ and the family $\{L_x:x\in \Pi\}$ consists of finite random point
processes on $\R^d$. Given $\Pi$, the family $\{L_x:x\in \Pi\}$ is i.i.d.\ with
the same distribution as a {\em typical cluster} $L_0$. We assume that
$$
\gamma:=E\card L_0
$$
is finite and positive, hence $\lambda=\lambda_\Pi\gamma$.

\begin{example}\label{exNS}\rm
Assume that
\begin{align}\label{cluster}
L_0=
\begin{cases}
\emptyset,&\quad\text{if $\eta=0$},\\
\{Y_{i,n}:i=1,\ldots,n\},& \quad \text{if $\eta=n\ge 1$},
\end{cases}
\end{align}
% Let
%\begin{align}\label{cluster}
%L_0=\left\{\begin{array}{ll} \emptyset &\quad {\rm if}\;\;\eta=0,\\
%\{Y_{i,n}:i=1,\ldots,n\}& \quad {\rm if}\;\;\eta=n\ge 1,
%\end{array}\right.
%\end{align}
where the random cluster size $\eta\ge 0$, and the random vectors
$Y_{i,n}$, $n\in\N$, $i=1,\ldots,n$, are independent.
Assume also that the $Y_{i,n}$  have the same (cluster point)
distribution $V$, say.
Then $N$ is called {\em Neyman-Scott process}. We always assume
that $\gamma=E[\eta]$ is positive and finite.
\end{example}

\begin{example}\label{exGP}\rm
Let $\eta$ be a $\{1,2\}$-valued random variable and assume that
\begin{align}\label{clusterGP}
L_0=
\begin{cases}
\{0\},&\quad\text{if $\eta=1$},\\
\{0,Y\},& \quad \text{if $\eta=2$},
\end{cases}
\end{align}
where $Y$ is a random vector independent of $\eta$.
Then $N$ is called \emph{Gauss-Poisson process}.
In this case a cluster  $L_0+x$ associated with a  parent
point $x$ say, contains $x$ and, with probability
$p:=P(\eta=2)$, also a secondary point $x$. Note that
the mean cluster size is given by $\gamma=1+p$.
\end{example}

In this paper we always assume that the {\em reduced second moment measure} 
\begin{align}\label{Lreduced}
E\Bigg[\sum_{\substack{x,y\in L_0\\x\ne y}}
\I\{x-y\in\cdot\}\Bigg]
\end{align}
of $L_0$ is absolutely continuous on $\R^d$. The well-known
second order properties of Poisson cluster point processes
(see e.g.\ \cite{DaVJ03})
easily imply that the measure defined at \eqref{reduced}
is absolutely continuous as well. For the Neyman-Scott process
of Example \ref{exNS} our assumption on $L_0$ is implied
by the absolute continuity of the cluster point distribution $V$.
For a Gauss-Poisson process it is sufficient to
assume that the seondary point $Y$ has an absolutely continuous
distribution.

In this section we fix a gauge body
$B$ that contains a non-empty neighborhood of $0$.
This way we exclude the trivial case $F_B(t)=0$,
$t\ge 0$. Let $\nu_B$ be the measure on $\R^d$ given by
\begin{align}\label{nuB}
\nu_B(\cdot):=d\int_{B^*}\I\{x/d_B(0,x)\in\cdot\}dx,
\end{align}
where $d_B(0,x):=d_B(0,\{x\})$. Using this measure
we can express the empty space hazard of $N$
as follows.

\begin{proposition}\label{prophazardPCB} The $B$-relative
empty space hazard of the Poisson cluster point process $Z=N$ is
given by
\begin{align}\label{hazardPC}
r_B(t,C)=
\lambda_\Pi t^{d-1}\int_{C} E\left[\int_{\R^d}
\I\{((L_0-x)\setminus \{0\})\cap (tu+tB)=\emptyset\}L_0(dx)\right]\nu_B(du).
\end{align}
\end{proposition}
{\em Proof.} Computing the left-hand side of \eqref{2.1}
for $K=\{0\}$ and $A=\{0\}\times C$ easily shows that
$$
db_dC_0(\{0\};B;\{0\}\times \cdot)=\nu_B(\cdot)
$$
and that $C_i(\{0\};B;\cdot)=0$ for $i\ge 1$. The result is then a
consequence of Proposition \ref{cclustehazard}  below.
\hfill\qed

\vspace{0.3cm}

If $L_0=\{0\}$ we have $N=\Pi$ and
\begin{equation}\label{hazardP}
r_B(t,C)=\lambda t^{d-1}\nu_{B}(C),
\end{equation}
in accordance with \eqref{densBoolean}.

If $B=B(1)$ then $\nu_B$ is the $(d-1)$-dimensional
Hausdorff measure on the unit sphere $S^{d-1}:=\partial B(1)$
and
\begin{align}\label{hazardPC3}
r_{B(1)}(t,C)= \lambda_\Pi t^{d-1}\int_{C}
E\left[\int_{\R^d}\I\{((L_0-x)\setminus \{0\})\cap
B(tu,t)=\emptyset\}L_0(dx)\right]\mathcal{H}^{d-1}(du).
\end{align}

\begin{example}\label{exNS2}\rm
Assume that $N$ is a Neyman-Scott process as defined
in Example \ref{exNS}. From \eqref{hazardPC} and
a straightforward calculation,
\begin{align}\label{NSHs}
r_B(t,C)=\lambda_\Pi t^{d-1}
\int_{\R^d}\int_{\R^d}g^\prime(P(Y_{1,1}-x\notin
tu+tB))V(dx)\nu_B(du),
\end{align}
where $g'$ is the derivative of the probability generating
function $g$ of $\eta$. This result generalizes formula (30)
in \cite{HLW02b}.
%This formula will be useful in the next
%proposition in order to compare Neyman-Scott processes.
\end{example}

\begin{example}\label{exGP2}\rm
Assume that $N$ is a Gauss-Poisson process as defined
in Example \ref{exGP2}. Then
\begin{align}\label{GPrB}
(\lambda_\Pi)^{-1}t^{1-d}r_B(t,C)=&(1-p)\nu_B(C)\\\notag
&+p\int_{C}P(Y\notin tu+tB) \nu_B(du)
+p\int_{C}P(-Y\notin tu+tB) \nu_B(du).
\end{align}
\end{example}

Let $N$ be a general Poisson cluster point process as above. It is
interesting to note that $t^{1-d}r_B(t,C)$ is monotone decreasing
in $t$. This is a direct consequence of \eqref{hazardPC}. Then it
is instructive to compare $r_B$ with the right side of
\eqref{hazardP}. From monotone convergence
\begin{align}\label{smallt}
\lim_{t\to 0}t^{1-d}r_B(t,C)=\lambda d\nu_B(C),
\end{align}
i.e.\ for small values of $t$ the empty space hazard of a Poisson
cluster point process behaves approximately like the empty space
hazard of a Poisson process with the same intensity.

Next we deal with the asymptotics of the empty space hazard for
large $t$. It turns out that it is the same as that of the Poisson
process $\Pi$ thinned at points $x$ where $L_x$ is empty, for
which  intensity equals $P(L_0\ne\emptyset)\lambda_\Pi$. This
means in a sense, that  points in clusters cannot be distinguished
from a very far distance, irrespective of any specific assumptions
on the typical cluster $L_0$. This is generalizing equation (22)
in \cite{HLW02b}. A weaker version of this
latter result has been rediscovered by Bordenave and Torrisi
\cite{BordTor07}.

\begin{proposition}\label{asymPC}
The $B$-relative empty space hazard of a Poisson cluster point process
satisfies
\begin{align}\label{larget}
\lim_{t\to \infty}t^{1-d}r_B(t,C)=P(L_0\ne\emptyset)\lambda_\Pi\nu_B(C).
\end{align}
\end{proposition}
{\em Proof.} The {\em tangential cone} (or {\em support cone})
$T(B,u)$ of $B$  at $u\in B$ is the closure
of $T'(B,u):=\{t(x-u):t>0, x\in B\}$, see Schneider \cite{Schneider1993}.
From \eqref{hazardPC} and monotone convergence
\begin{align}\label{larget2}
\lim_{t\to \infty}t^{1-d}r_B(t,C)=\lambda_\Pi\int_{C}
E\left[\int\I\{((L_0-x)\setminus \{0\})\cap
T'(B,-u)=\emptyset\}L_0(dx)\right]\nu_B(du).
\end{align}
Since the measure \eqref{Lreduced} is absolutely continuous
it follows that $Y_0$ a.s.\ does not intersect the
boundary of  $T(B,-u)$. Hence
\begin{align}\label{target2}
\lim_{t\to \infty}t^{1-d}r_B(t,C)=\lambda_\Pi\int_{C}
E\left[\int\I\{((L_0-x)\setminus \{0\})\cap
T(B,-u)=\emptyset\}L_0(dx)\right]\nu_B(du).
\end{align}
Let us now fix for a moment a {\em regular}
boundary point $u$ of $B^*$. This means that
$B$ has a unique supporting hyperplane at $-u$,
see Schneider \cite{Schneider1993}.
Then $T(B,-u)$ is the supporting half-space of $B$ at $-u$,
see Section 2.2 in Schneider \cite{Schneider1993}.
It now follows as in \cite[Section 6.5]{HLW02b} that
$$
E\left[\int\I\{((L_0-x)\setminus \{0\})\cap
T(B,-u)=\emptyset\}L_0(dx)\right]=P(L_0\ne\emptyset).
$$
It remains to note that $\nu_B$-a.a.\ $u\in\R^d$ are regular
boundary points of $B^*$. This follows from the fact that $\nu_B$
is absolutely continuous with respect to $(d-1)$-dimensional
Hausdorff measure on $\partial B^*$ and Theorem 2.2.4 in Schneider
\cite{Schneider1993}. 
\hfill\qed

\vspace*{0.3cm}

%Can we give a formula for Hawkes point
%processes in terms of the usual parameters of Hawkes?})

%with random cluster sizes $\eta$ and $\tilde\eta$, respectively,
%but with the same distribution $V$ of the random vectors
%$Y_{i,n}$. In order to compare variability of $\eta$, and
%$\tilde\eta$ typically in use are the following orderings:
%$\eta<_{icv} \tilde\eta$ iff $E[\psi(\eta)]\le
%E[\psi(\tilde\eta)]$ for all non-decreasing concave functions
%$\psi$. If $E(\eta)= E(\tilde\eta)$ then $\eta<_{icv} \tilde\eta$
%iff $\tilde\eta<_{cx} \eta$, for the convex ordering $<_{cx}$
%defined analogously as $<_{icv}$, but using convex functions
%$\psi$. The following proposition formalizes intuition that more
%variable number of secondary points in clusters may result in
%larger  empty space in terms of the hazard rate ordering.
%Intuitively, more variability results in more clustered patterns,
%which results in more empty space. This scenario works for
%positive variables $\eta,\tilde\eta$, in which case thinning is
%not present.

We shall now consider two Poisson cluster processes
$N$ and $\tilde N$ having the same intensity,
satisfying the technical assumption
formulated at \eqref{Lreduced} and with
$B$-relative hazard rates
$r_B$ and $\tilde r_B$, respectively. We will establish
several sets of assumptions implying the hazard rate
ordering \eqref{2.8}.
In all cases $\tilde N$ will have more clustered patterns
in a sense. Therefore one
may expect more empty space in the Poisson cluster
process $N$. Still it is somewhat surprising that
this happens in the strong sense of \eqref{monhazard}.
In Section \ref{secpclgg} we will state the corresponding results
for Poisson cluster grain models.
Our first ordering result, dealing with Neyman Scott processes,
is a special case of
Proposition \ref{propNS} which will be proved later.

Let $\eta$ and $\tilde\eta$
be two counting variables (i.e. taking values in $\{0,1,2,\ldots\}$) and let $\eta_l$
be the (shifted) {\em length-biased} version of $\eta$. This means
that $\eta_l$ has distribution
$E[\eta]^{-1}E[\eta\I\{\eta-1\in\cdot\}]$.
Denoting by $\tilde\eta_l$ the length-biased version of $\tilde\eta$,
we define the \emph{length biased probability generating functions ordering }
\begin{align}\label{biasedordering}
\eta<_{l-g}\tilde\eta
\end{align}
by
\begin{align*}
E[s^{\eta_l}]
&\ge E[s^{\tilde\eta_l}],\quad s\in[0,1].
\end{align*}
This means that $\eta_l$ is smaller than $\tilde\eta_l$
in the the generating function order (see Stoyan \cite{St83}, Section 1.8.)
Note that $\eta<_{l-g}\tilde\eta$ is equivalent to 
$$
E[\eta]^{-1}E[\I\{\eta\ge 1\}\eta a^{\eta-1}]
\ge E[\tilde\eta]^{-1}E[\I\{\tilde\eta\ge 1\}\tilde\eta a^{\tilde\eta-1}],\quad a\in[0,1].
$$  
Another way of expressing this relation is by the generating functions 
of $\eta$ and $\tilde\eta$.   Denote by
$$
g_\eta(s):=E[s^\eta],\quad s\in[0,1],
$$
the probability generating function of $\eta$.
Then  \eqref{biasedordering} means that
\begin{align}\label{gfbiasedordering}
E[\eta]^{-1}g'_\eta(s)\ge E[\tilde\eta]^{-1}g'_{\tilde \eta}(s),
\quad s\in[0,1].
\end{align}
The relation \eqref{biasedordering} does not imply that the corresponding 
expected values $E[\eta]$,$E[\tilde\eta]$ are ordered.  
For example if $\eta\equiv 1$,  $\eta_1\equiv 2$, and $\eta_2$ equals 0 
with probability $1/2$ and $1$ with probability $1/2$ 
then $\eta<_{l-g}\eta_1$, and $\eta<_{l-g}\eta_2$, but $E[\eta]=1<E[\eta_1]=2$, 
and $E[\eta]=1>E[\eta_2]=1/2.$ However, if  $E[\eta]=E[\tilde\eta]$ 
then \eqref{biasedordering}  implies that $\Var[\eta]\le \Var[\tilde\eta]$, 
therefore this relation is a variability ordering.

\begin{proposition}\label{propNSpp}
Assume that $N$ and $\tilde{N}$ are Neyman-Scott processes
with cluster sizes $\eta$ and $\tilde\eta$, respectively,
having the same distribution $V$ of cluster points.
If $N$ and $\tilde{N}$ have the same intensity and
\begin{align*}
\eta<_{l-g}\tilde\eta
\end{align*}
then
\begin{equation}\label{monhazard}
Z<_{h-B}\tilde Z
\end{equation}.
\end{proposition}

\begin{example}\label{ex5}\rm Let $N$ and $\tilde{N}$ be as
in Proposition \ref{propNSpp} and assume that
$\eta\equiv 1$. That is  $N$ is again a Poisson process with intensity
$1$ since the  points of the original Poisson process $\Pi$
are independently shifted. Let $a\in[0,1]$. The length 
biased variable $\eta_l\equiv 0$, and its generating function $E[s^{\eta_l}]\equiv 1$, 
therefore for each $\tilde \eta$ we have $1\equiv\eta<_{l-g}\tilde\eta$, 
and $Z<_{h-B}\tilde Z$.
This example shows that within the class of Neyman-Scott processes with 
fixed intensity the stochastically smallest empty space appears for pure 
Poisson germ processes.
\end{example}

\begin{example}\label{ex6}\rm Let $N$ and $\tilde{N}$ be as
in Proposition \ref{propNSpp} and assume that
$\eta$ and $\tilde\eta$ are Poisson distributed
with parameter $c$ and $\tilde c$, respectively.
(Then the clusters are finite Poisson processes.)
We have for any $a\in[0,1]$ that
$$
E[\eta]^{-1}E[\I\{\eta\ge 1\}\eta a^{\eta-1}]=e^{-c(1-a)}.
$$
Therefore \eqref{biasedordering} (and hence $Z<_{h-B}\tilde Z$) holds iff
$c\le \tilde c$.  Note that  for Poisson distributed $\eta$ the length-biased
$\eta_l$ have the same (Poisson) distribution. Proportionally more Poisson cluster 
points in a Neyman-Scott process indeed generate more clustering and 
lead to a stochastically larger empty space.
\end{example}

\begin{example}\label{ex7}\rm  Assume that
$\eta$ and $\tilde\eta$ are binomial distributed
with parameter $(n,p)$ and $(\tilde n,\tilde p)$, respectively.
Since $g_\eta(s)=((1-p)+ps)^n$, it follows that
\eqref{biasedordering} is equivalent to
\begin{align}
((1-p)+ps)^{n-1}\ge ((1-\tilde p)+\tilde ps)^{\tilde n-1},\quad s\in[0,1].
\end{align}
If, for instance $n=\tilde n$, then this inequality is
implied by $p\le \tilde p$.
\end{example}
\begin{example}\label{ex8}\rm
 Assume that
$\eta$ and $\tilde\eta$ are negative binomial distributed
with parameter $(p,r)$ and $(\tilde p,\tilde r)$, respectively. The corresponding 
length-biased variables $\eta_l$ , $\tilde\eta_l$ are again negative binomial 
distributed with parameters $(p,r+1)$ and $(\tilde p,\tilde r+1)$, respectively. 
For $p=\tilde p$, if $r\le \tilde r$ then $\eta<_{l-g}\tilde\eta$, 
and for $r=\tilde r$, if $p\ge \tilde p$ then $\eta <_{l-g}\tilde\eta$. 
This is a special case of a more general setting. 
If $\eta =\sum_{i=1}^{\kappa } \vartheta_i$, 
and $\tilde\eta =\sum_{i=1}^{\tilde\kappa} \tilde\vartheta_i$, 
for iid variables $\{\vartheta_i\}_{i\ge 1}$, and independent of $\kappa$ 
(all variables taking on natural values) then 
$\vartheta_i <_{l-g}\tilde \vartheta_i$, and $\kappa <_{l-g}\tilde\kappa$ implies  
$\eta <_{l-g}\tilde\eta$.
\end{example}

The next result is a special case of Proposition \ref{propGP}.

\begin{proposition}\label{propGPpp}
Assume that $N$ and $\tilde{N}$ are Gauss-Poisson processes
having the same distribution $V$ of the secondary point.
Assume that $N$ and $\tilde{N}$ have the same intensity and
that the probabilities $p$ and $\tilde p$ for having
a secondary point satisfy  $p\le \tilde p$. Then
$Z<_{h-B}\tilde Z$.
\end{proposition}

In our next result we will multiply each point of the typical
cluster $L_0$ of $N$ with a random
variable $W\le 1$, i.e.\ we assume that the typical
cluster $\tilde L_0$ of $\tilde N$ has the
distribution of $WL_0=\{Wx:x\in L_0\}$. Compared with $\tilde L_0$, the
points of $L_0$ are more spread out.
Note that we allow any sort of dependence between $L_0$ and $D$.

\begin{proposition} \label{spreadoutpp}
Consider two Poisson cluster processes $N$ and $\tilde N$ based
on the same Poisson process $\Pi$ and typical clusters
$L_0$ and $\tilde L_0$, respectively.
Assume that $\tilde L_0$ is distributed as $WL_0$ for some random variable
$W\le 1$.  Then $Z<_{h-B}\tilde Z$.
\end{proposition}
{\em Proof.} Take $x_1,\ldots,x_n\in \R^d$ and
set $\psi:=\{x_j:j=1,\ldots,n\}$.
Let $w\le 1$ and define $\tilde \psi:=\{wx_j:j=1,\ldots,n\}$.
Let $u\in\partial B^*$ and $t\ge 0$.
In view of \eqref{hazardPC} it is sufficient
to show that
\begin{align}\label{8}
\I\{((\tilde\psi-wx_i)\setminus \{0\})\cap (tu+tB)=\emptyset\}
\le \I\{((\psi-x_i)\setminus \{0\})\cap (tu+tB)=\emptyset\}
\end{align}
holds for any $i\in\{1,\ldots,n\}$. Assume that the
left-hand side of \eqref{8} equals $1$. This is equivalent
to $w(x_j-x_i)\notin tu+tB$ for all $j\ne i$. Since
$u\in\partial B^*$ we have that $0\in\partial(tu+tB)$.
Therefore the convexity of $tu+tB$ and $w\le 1$ imply that also
$x_j-x_i\notin tu+tB$ for all $j\ne i$.
This shows \eqref{8} and hence the proposition. 
\hfill\qed

\vspace{0.3cm}
Recall that a random variable
$\tilde{W}$ is {\em stochastically smaller}
than another r.v.\ $W$, if
$P(\tilde W>t)\le P(W>t)$ for all $t\ge 0$.

\begin{proposition} \label{soscaling}
Let $N$ and $\tilde N$ be as in Proposition \ref{spreadoutpp}.
Let $L$ be a point process and assume that $L_0$ is
distributed as $WL$ for a random variable $W>0$ independent
of $L$. Assume that $\tilde{L}_0$ is
distributed as $\tilde{W}L$ for a random variable $\tilde{W}\ge 0$
independent of $L$. If $\tilde{W}$ is stochastically smaller
than $W$ then $Z<_{h-B}\tilde Z$.
\end{proposition}
{\em Proof.} By inverse coupling based on an uniformly
distributed random variable that it independent of
$L$ we can assume that $(W,\tilde W)$ is independent
of $L$ and that $\tilde W\le W$ everywhere on the underlying
probability space. Since $\tilde{L}_0=\tilde{W}/W L_0$,
we can apply Proposition \ref{spreadoutpp} with
$W$ replaced by $\tilde{W}/W\le 1$.
\hfill\qed

\section{Results for Poisson cluster germ-grain models}\label{secpclgg}

We now return to a general case of germ-grain models where $N$ is
a Poisson cluster point processes as in the previous section, and
grains $X_n$, located at all points of $N$, are compact and
convex. In the next two sections we fix a structuring element
$B$ such that our general assumption made after
\eqref{reduced} holds, that is $X_1$
and $B^*$ are a.s.\ in general relative position. Alternatively, 
we shall treat such a process as a Boolean
model with non-convex grains which are
$$
\bigcup_{y\in L_x}(X(y)+y), \quad x\in \Pi.
$$

For a given finite point pattern
$\psi=\{x_n:n=1,\ldots,m\}$, denote by $\Gamma(\psi,\cdot)$ the distribution of
the random closed set $\cup_{n=1}^m (X_n+x_n)$, where $X_1,\ldots,X_m$ are
independent with distribution of $X_0$. We also set
$\Gamma(\emptyset,\cdot):=\delta_{\emptyset}$. Let $Y_0$ be a random closed set
with distribution
\begin{align}\label{mu}
\mu(\cdot):=\frac{1}{\gamma}E\left[\sum_{y\in L_0}
\int\I\{A\in \cdot\}\Gamma((L_0-y)\setminus\{0\},dA)\right].
\end{align}
This probability measure is describing the distribution of the
germ-grain model associated with a typical cluster as seen from a
randomly chosen cluster point, after removing the grain around the
chosen point. We assume that $Y_0$ and the typical grain $X_0$ are
independent and define for $i\in\{0,\ldots,d-1\}$, $t\ge 0$, and a
Borel set $C\subset\R^d$
\begin{align}\label{HiB}
K_{i,B}(t,C):=\frac{1}{\bar V_{i,B}}
E\left[\int\I\{d_B(x+tu,Y_0)>t\}\I\{u\in C\}C_{i}(X_0;B;d(x,u))\right].
\end{align}
The following proposition yields again an explicit formula for the empty space
hazard of a Poisson cluster germ-grain model, this time in terms of $X_0$ and
$Y_0$ describing locally a cluster.

\begin{proposition}\label{cclustehazard}
The $B$-relative empty space hazard of a Poisson cluster
germ-grain model with compact, convex grains is given by
\begin{align}\label{densclust}
r_B(t,C)=\sum^{d-1}_{i=0}(d-i)t^{d-i-1}b_{d-i}\lambda\bar{V}_{i,B}
K_{i,B}(t,C).
\end{align}
%where $K_{i,B}$ is defined by \eqref{HiB}.
\end{proposition}
{\em Proof.} Our aim is to use \eqref{dens1}. To do so we
recall that the Palm probability measure $P^0_N$
of a Poisson cluster process $N$ satisfies
\begin{align}\label{ff}
P^0_N(N\in\cdot)= E\left[\int\I\{N\cup\psi\in\cdot\}Q^0_{L_0}(d\psi)\right],
\end{align}
where
$$
Q^0_{L_0}(\cdot):=\gamma^{-1}E\left[\sum_{y\in L_0} \I\{L_0-y\in\cdot\}\right]
$$
is the {\em Palm distribution} of the  typical cluster $L_0$, see
e.g.\ Stoyan, Kendall and Mecke \cite{SKM95}. As in  Last and
Holtmann \cite{LH99}, this implies that
\begin{align*}
P^0_N (d_B(x+tu,Z^{!})>t)
=(1-F_B(t))\iint\I\{d_B(x+tu,A)>t\}\Gamma(\psi\setminus\{0\},dA)Q^0_{L_0}(d\psi).
\end{align*}
By definition \eqref{GiB} and definition of $Y_0$ this means that
\begin{align*}
\frac{G_{i,B}(\infty,C)-G_{i,B}(t,C)}{1-F_B(t)}
=\frac{1}{\bar V_{i,B}}E\left[\int\I\{d_B(x+tu,Y_0)>t\}\I\{u\in C\}
\bar C_{i,B}(d(x,u))\right].
\end{align*}
By definition  \eqref{HiB}, the above right-hand side equals
$K_ {i,B}(t,C)$. Inserting this into \eqref{dens1} gives
the asserted equation \eqref{densclust}. 
\hfill\qed

%\subsection{Asymptotic properties of empty space hazard}

\vspace*{0.3cm}

Our next proposition deals with the asymptotic behavior of
$r_B(t,C)$ as $t\to 0$ or $t\to\infty$, respectively. Note that
$t^{1-d}r_B(t,C)$ is monotone decreasing in $t$. This is a direct
consequence of \eqref{HiB}. Recall the definition of the
tangential cone given in the proof of Proposition \ref{asymPC}.

\begin{proposition}\label{casymptotics}
The $B$-relative empty space hazard $r_B(t,C)$ of a Poisson
cluster germ-grain model with compact, convex grains satisfies
\begin{align}\label{as1}
\lim_{t\to 0}r_B(t,C)&=2\lambda_\Pi E\left[\int_{\R^d\times C}
\I\{x\notin Y_0\}C_{d-1}(X_0;B;d(x,u))\right]\\ \label{as2}
\lim_{t\to \infty}t^{1-d}r_B(t,C)
&=db_d\lambda_\Pi E\left[\int_{\R^d\times C}\I\{(x+T(B,-u))\cap
Y_0=\emptyset\}C_{0}(X_0;B;d(x,u))\right].
\end{align}
\end{proposition}
{\em Proof.} The proof is similar to the proof of \eqref{smallt}
and \eqref{target2}. \hfill\qed

\begin{remark}\label{positive}\rm One might wonder whether
the right-hand side of \eqref{as2} is positive for
$C=\R^d$. A simple sufficient
condition is that $P(\card L_0=1)>0$,
because in this case $Y_0$ is with positive probability empty.
Another sufficient condition is to assume that $B$ is
smooth (any boundary point has a unique supporting hyperplane)
and that the diameter of the typical grain can take arbitrary small
positive values with positive probability. These assumptions
would allow to apply the method of \cite[Section 6.5]{HLW02b}
on a set of positive probability. We do not go into further
details.
\end{remark}

\begin{remark}\label{ligh-tailed}\rm
The second assertion of the previous proposition shows in
particular that $F_B$ is {\em light-tailed}, i.e.\ has a finite
exponential moment.
\end{remark}

%\subsection{Comparison of Poisson cluster germ-grain models}

Now we shall compare the $B$-relative empty space hazard $r_B$ of
$Z$ with that of another Poisson cluster germ-grain model
$\tilde{Z}$ with the same intensity $\lambda$ of germs and the
same typical grain $X_0$. Both underlying germ processes are
assumed to satisfy the technical assumption formulated at \eqref{Lreduced}.
We denote the characteristics of $\tilde{Z}$ by $\tilde{\Pi}$,
$\tilde{L}_0$, $\tilde{K}_{i,B}$, $\tilde{r}_B$ etc. We begin with
a direct consequence of Theorem \ref{cclustehazard}.

\begin{proposition}\label{general}
Let $C\in\mathcal{B}^d$ and $t\ge 0$ such that
\begin{align}\label{monass}
K_{i,B}(t,C)\ge \tilde K_{i,B}(t,C).
\end{align}
Then \eqref{monhazard8} holds.
\end{proposition}

\vspace{0.3cm}
An immediate consequence of the above proposition
is that the relative empty space hazard of a Boolean model is
always greater than that of a Poisson cluster germ-grain model
having the same germ intensity and the same typical grain:

\begin{corollary}\label{cBooleanempty}
Assume that $Z$ is a Boolean model with typical convex, compact
grains distributed as $X_0$, and $\tilde Z$ a Poisson cluster
germ-grain model with equal intensity, and also typical grains
distributed as  $X_0$. Then $Z<_{h-B}\tilde Z$.
\end{corollary}
{\em Proof.} Since $Z$ is a Boolean model, we have
$L_0=\{0\}$ and $\gamma=1$. Hence, if $\bar{V}_{i,B}>0$,
\begin{align}\tag*{\qed}
K_{i,B}(t,C)=\bar{V}_{i,B}^{-1}\bar C_{i,B}(\R^d\times C)
\ge \tilde K_{i,B}(t,C).
\end{align}

\vspace{0.3cm}
While a Boolean model has stochastically smaller empty space
than a related Poisson cluster germ-grain model,
it has a greater volume fraction. 
Under a different set of assumptions (more specific Poisson-cluster
processes and deterministic but possibly non-convex grains)
the result was proved in Section 3.8 of \cite{Hall88}.

\begin{proposition}\label{probvolumfraction}
Under the assumptions of Corollary \ref{cBooleanempty}
\begin{align}\label{4.6}
P(0\in Z)\ge P(0\in\tilde Z).
\end{align}
\end{proposition}
{\em Proof.} The volume fraction of $Z$ is given by
\eqref{vfBoolean}.
On the other hand, $\tilde Z$ is also a Boolean model, but based
on the Poisson process $\tilde{\Pi}$ and with typical (possibly
non-convex) grain
$$
\tilde X_0:=\bigcup_{x\in \tilde{L}_0}X(x)+x,
$$
where $\tilde{L}_0$ is the typical cluster associated with $\tilde
Z$ and, given $\tilde{L}_0$, the family $\{X(x):x\in
\tilde{L}_0\}$ consists of independent random closed sets with the
same distribution as $X_0$. Therefore we obtain from
\eqref{vfBoolean}
\begin{align}\label{4.8}
P(0\in \tilde Z)=1-\exp[-\lambda_{\tilde\Pi} E V_d(\tilde X_0)].
\end{align}
We have
\begin{align*}
E V_d(\tilde X_0)\le E \sum_{x\in \tilde L_0}V_d(X(x)) =E
\sum_{x\in \tilde L_0}E[V_d(X(x))|\tilde L_0]=\gamma \bar V_d.
\end{align*}
Inserting this into \eqref{4.8} and comparing with
\eqref{vfBoolean}, yields the assertion.
\hfill\qed

%\vspace{0.3cm}
\begin{remark}\label{rem1}\rm Consider the hypothesis of
Proposition \ref{probvolumfraction}. The proof of this proposition
shows that we have equality in \eqref{4.6} iff
\begin{align}\label{4.11}
V_d(X(x)\cap X(y))=0\quad x,y\in \tilde L_0,x\ne y\quad
P-\text{a.s.}
\end{align}
This is, for instance, the case if the cluster points have minimal
distance $2t_0$  from each other for some $t_0>0$, and $X_0$ is
a.s.\ contained in the ball $B(t_0)$.
\end{remark}

Next we prove a more general germ-grain version
of Proposition \ref{propNSpp}.

\begin{proposition}\label{propNS}
Consider two germ-grain models with the same typical grains $X_0$
such  that $N$ and $\tilde{N}$ are Neyman-Scott processes with
cluster sizes $\eta$ and $\tilde\eta$, respectively, and the same
cluster point distribution $V$. If $N$ and $\tilde{N}$ have the same intensity and
\begin{align}
\eta<_{l-g}\tilde\eta
\end{align}
then $Z<_{h-B}\tilde Z$.
\end{proposition}
{\em Proof.} By Proposition \ref{general} it suffices to show
\begin{align}\label{suff}
P(d_B(x+tu,Y_0)>t)\ge P(d_B(x+tu,\tilde Y_0)>t)
\end{align}
for all $x,u\in\R^d$ and all $t\ge 0$. Letting $B':=tB+x+tu$
and using the definition \eqref{mu} of the distribution of $Y_0$, we
obtain that
\begin{align*}
E[\eta]P(Y_0\cap B'=\emptyset)
&=E\left[\sum_{y\in L_0}
\I\{A\cap B'=\emptyset\}\Gamma((L_0-y)\setminus\{0\},dA)\right]\\
&=P(\eta=1)+\sum^\infty_{n=2}P(\eta=n)n\int f(y)^{n-1}V(dy),
\end{align*}
where $f(y):=\int P((X_0+z-y)\cap B'=\emptyset)V(dz)$ and where the
second identity comes from a straightforward calculation using
the definition of the typical cluster $L_0$ of a Neyman-Scott process.
By Fubini's theorem this means that
\begin{align*}
E[\eta]P(Y_0\cap B'=\emptyset)=\int E[\eta f(y)^{\eta-1}]V(dy).
\end{align*}
We can now use our assumption \eqref{biasedordering} to derive
\begin{align*}
P(Y_0\cap B'=\emptyset)\ge
E[\eta]^{-1}\int E[\tilde\eta f(y)^{\tilde\eta-1}]V(dy).
\end{align*}
Reversing the above steps,
we get \eqref{suff} and hence the asserted result.
\hfill\qed

Our next result generalizes Theorem 5.4 in \cite{LH99}.

\begin{proposition}\label{propGP}
Consider two germ-grain models with the same typical grains $X_0$
based on Gauss-Poisson processes $N$ and $\tilde{N}$.
Assume that $N$ and $\tilde{N}$ have the same intensity and
that the probabilities $p$ and $\tilde p$ for having
a secondary point satisfy  $p\le \tilde p$. Then
$Z<_{h-B}\tilde Z$.
\end{proposition}
{\em Proof.} Fix $t\ge 0$ and $C\in\mathcal{B}^d$.
From the defining properties of a Gauss-Poisson process
(see Example \ref{exGP2}) we have
\begin{align*}
\bar{V}_{i,B}K_{i,B}(t,C)=&
\frac{1-p}{1+p}\bar C_{i,B}(\R^d\times C)
+\frac{p}{1+p}\int_{\R^d\times C} a(x,u)\bar C_{i,B}(d(x,u))\\
&+\frac{p}{1+p}\int_{\R^d\times C} b(x,u)\bar C_{i,B}(d(x,u)),
\end{align*}
where
\begin{align*}
a(x,u):=P(d_B(x+tu,Y+X_0)>t),\qquad b(x,u):=P(d_B(x+tu,-Y+X_0)>t),
\end{align*}
where $Y$ and $X_0$ are independent.
Therefore, by Proposition \ref{general}, it suffices to show
that
\begin{align*}
\frac{1-p}{1+p}+\frac{pa}{1+p}+\frac{pb}{1+p}
\ge \frac{1-\tilde p}{1+\tilde p}+
\frac{\tilde pa}{1+\tilde p}+\frac{\tilde pb}{1+\tilde p},
\end{align*}
for all $a,b\in[0,1]$. Simple algebra shows that this
inequality is equivalent to
$$
2\tilde p-a\tilde p-b\tilde p\ge 2p-ap-bp.
$$
The latter is implied by our assumption $p\le \tilde p$.
\hfill\qed

\section{Results for Mixed Poisson germ-grain models}

In this section we consider a germ-grain model $Z$ based
on a {\em mixed Poisson process} $N$. This means that
there is a random variable $\Lambda\ge 0$ such that
the conditional distribution of $N$ given $\Lambda$
is that of a stationary Poisson process with intensity $\Lambda$.
We assume that  $E[\Lambda]$ (the
            intensity of $N$)
            is positive and finite.

It is convenient to use the notation
\begin{align*}
H_B(t)=E[V_d(X_0+tB^*)],\quad t\ge 0.
\end{align*}

\begin{proposition}\label{cmixedhazard}
The $B$-relative empty space hazard of a mixed Poisson
germ-grain model with compact, convex grains is given by
\begin{align}\label{densmixed}
r_B(t,C)=\sum^{d-1}_{i=0}(d-i)t^{d-i-1}b_{d-i}
E[\exp[-\Lambda H_B(t)]]^{-1}E[\Lambda\exp[-\Lambda H_B(t)]]
\bar{S}_{i,B}(t,C).
\end{align}
\end{proposition}
{\em Proof.} Again we will use \eqref{dens1}. To do so, we
note that the Palm probability measure $P^0_N$
of a mixed Poisson process $N$ satisfies
\begin{align}\label{Palmmixed}
P^0_N((\Lambda,N)\in\cdot)=E[\Lambda]^{-1}E[\Lambda\I\{(\Lambda,N\cup\{0\})\in\cdot\}].
\end{align}
This  formula can be derived by conditioning and using
the properties of a Poisson process.
Since, moreover, the conditional distribution
$P^0_N(Z^!\in\cdot|\Lambda)$ (cf.\ \eqref{Z!} for the
definition of the random set $Z^!$) is that
of a Boolean model with germ intensity $\Lambda$, we
obtain for all $(x,u)\in\R^d\times\R^d$ that
\begin{align*}
P^0_N (d_B(x+tu,Z^{!})>t)&=
E[\Lambda]^{-1}E[\Lambda P(d_B(x+tu,Z)>t\mid\Lambda)]\\
&=E[\Lambda]^{-1}E[\Lambda P(d_B(0,Z)>t\mid\Lambda)]\\
&=E[\Lambda]^{-1}E[\Lambda \exp[-\Lambda H_B(t)]],
\end{align*}
where we have used \eqref{CFBoolean} to obtain
the last identity. Again by conditioning and
\eqref{CFBoolean} we have that $1-F_B(t)=E[\exp[-\Lambda H_B(t)]]$.
Inserting our findings into the general formula \eqref{dens1}
yields the assertion \eqref{densmixed}.
\hfill\qed 

\vspace*{0.3cm}

In order to
state some stochastic ordering consequences of Proposition \ref{cmixedhazard}
we introduce a stochastic order using  cumulants.
For two nonnegative random variables $\Lambda$, $\tilde\Lambda$ we say that they are ordered in the first cumulant order and write $\Lambda<_{cum}\tilde\Lambda$ if for the corresponding cumulant generating functions $C_{\Lambda}$,  $C_{\tilde\Lambda}$, taking the first derivatives we have  $C'_{\Lambda}(s)\ge C'_{\tilde\Lambda}(s), \ s\in [-1,0].$
Note that $\Lambda<_{cum}\tilde\Lambda$ is equivalent  to
\begin{align}\label{4}
E[\exp[-\Lambda s]]^{-1}E[\Lambda\exp[-\Lambda s]]
\ge E[\exp[-\tilde\Lambda s]]^{-1}E[\tilde\Lambda\exp[-\tilde\Lambda s]],\quad s\ge 0.\end{align}   The left-hand side of \eqref{4}
is the logarithmic derivative of the
Laplace transform $s\mapsto E[\exp[-\Lambda s]]$. It is also
the hazard rate of the distribution function $G_\Lambda$, defined
by
\begin{align*}
G_\Lambda(s):=1-E[\exp[-\Lambda s]],\quad s\ge 0.
\end{align*}
This is a mixture of exponential distributions. Equation \eqref{4}
then means that the corresponding variables are ordered
in the hazard rate order, i.e.  $G_\Lambda<_{h}G_{\tilde\Lambda}.$
Note that for $\Lambda$, $\tilde\Lambda$ with equal expected values  $\Lambda<_{cum}\tilde\Lambda$ implies that $\Var[\Lambda]\le \Var[\tilde\Lambda]$, therefore, similarly to the relation $<_{l-g}$, the relation $<_{cum}$ is a variability ordering in the case of a fixed expectation.

An immediate consequence of Proposition \ref{cmixedhazard}
is the following counterpart of Propositions \ref{propNS}
and \ref{propGP}. We use similar notation.
Again, intuitively speaking, more variability in the 
mixed Poisson model results in a stochastically greater empty space 
(a stochastically larger clustering).

\begin{proposition}\label{propmixed}
Consider two germ-grain models with the same typical grains $X_0$
based on mixed Poisson processes  $N$ and $\tilde{N}$ with
random intensities $\Lambda$ and $\tilde\Lambda$, respectively.
Assume that
$\Lambda<_{cum}\tilde\Lambda$ then $Z<_{h-B}\tilde Z$.
\end{proposition}

\begin{example}\label{exGamma}\rm Let $N$ and $\tilde{N}$ be as
in Proposition \ref{cmixedhazard} and assume that
$\Lambda$ is Gamma distributed with shape and scale parameter
$\alpha>0$ and $\beta>0$, respectively. This means that
$\Lambda$ has density
$\beta^\alpha \Gamma(\alpha)^{-1}x^{\alpha-1}\exp[-\beta x]$.
The Laplace transform of $\Lambda$ can be computed as
\begin{align*}
E[\exp[-s\Lambda]]=\frac{\beta^\alpha}{(\beta+s)^{\alpha}},
\end{align*}
while an equally easy calculation gives
\begin{align*}
E[\Lambda\exp[-s\Lambda]]=\frac{\alpha}{\beta+s}
\frac{\beta^\alpha}{(\beta+s)^{\alpha}}.
\end{align*}
Assume now that $\tilde\Lambda$ is Gamma distributed with parameters
$\tilde\alpha$ and $\tilde\beta$, respectively. Then assumption \eqref{4}
means that $\alpha/(\beta+s)\ge \tilde\alpha/(\tilde\beta+s)$
holds for all $s\ge 0$. This is equivalent to
\begin{align}\label{34}
\alpha\ge \tilde\alpha,\quad
\frac{\alpha}{\beta}\ge \frac{\tilde\alpha}{\tilde \beta}.
\end{align}
Depending on whether or not $\beta\ge \tilde\beta$, only
one of these equations is relevant.
By Proposition \ref{propmixed}, \eqref{34} implies the empty space hazard 
ordering \eqref{monhazard}.
Assume for instance that $\tilde\Lambda$ is exponentially
distributed with mean $1$, i.e.\ $\tilde\alpha=\tilde\beta=1$
and assume furthermore that $\Lambda$ has also mean $1$, that is
$\alpha=\beta$. Then \eqref{34} is equivalent to
$\alpha\ge \tilde\alpha$. Note that the variance of $\Lambda$
satisfies  $\Var[\Lambda]=\alpha/\beta^2=1/\beta\le \Var[\tilde\Lambda]=1$ if
$\alpha\ge 1$.
\end{example}

As in the Poisson cluster case it follows
that the relative empty space hazard of a Boolean model is
greater than that of a mixed Poisson germ-grain model
with the same germ intensity.

\begin{corollary}\label{cBooleanempty2}
Assume that $Z$ is a Boolean model with typical convex, compact
grains distributed as $X_0$, and $\tilde Z$ a mixed Poisson
germ-grain model with equal intensity, and also typical grains
distributed as  $X_0$. Then $Z<_{h-B}\tilde Z$.
\end{corollary}
{\em Proof.} Let $\lambda$ denote the germ intensity
of the Boolean model $Z$ and let $\tilde\Lambda$ be the
random intensity of the mixed Poisson process underlying
$\tilde Z$. It is assumed that $E[\tilde \Lambda]=\lambda$.
We check that condition \eqref{4} holds with $\Lambda\equiv\lambda$.
This condition means that
\begin{align*}
\lambda E[\exp[-\tilde\Lambda s]]
\ge E[\tilde\Lambda\exp[-\tilde\Lambda s]],\quad s\ge 0.
\end{align*}
In other words: the covariance between $\tilde\Lambda$
and $-\exp[-\tilde\Lambda s]$ has to be non-negative.
This fact follows from a very well-known statement that a single 
random variable is associated, see Esary et al. \cite{Esary67}.
\hfill\qed

\medskip
For completeness we provide the mixed Poisson analogue of
Corollary \ref{probvolumfraction}. The result can be found in
Section 3.8 of \cite{Hall88}
for the more general case of stationary Cox processes
with an absolutely continuos intensity measure.
Although our proof below can be extended to arbitrary stationary
Cox processes we stick to the mixed Poisson case.

\begin{proposition}\label{probvfracmixed}
Under the assumptions of Corollary \ref{cBooleanempty2}
\begin{align}\label{4.23}
P(0\in Z)\ge P(0\in\tilde Z).
\end{align}
\end{proposition}
{\em Proof.} By conditioning and \eqref{vfBoolean},
\begin{align*}
1-P(0\in \tilde Z)=E[\exp[-\Lambda E[V_d(X_0)]]].
\end{align*}
By Jensen's inequality this is bounded from below
by $\exp[-\lambda E[V_d(X_0)]]$. This lower bound
is just $1-P(0\in  Z)$.
\hfill\qed

\section{Concluding remarks}\label{problems}

We have derived several variability properties
            of the empty space function of Poisson cluster and
            mixed Poisson germ-grain models.
            It would be worthwhile to study also other classes
            of germ processes.
            Another interesting task is to
            find a good notion of {\em spread out} for a
            finite point process (with respect to the origin).
            Proposition \ref{propGPpp} and Proposition \ref{spreadoutpp}
            should be both special cases of the same principle.
            The first proposition is generalized by Proposition \ref{propGP}.
            We believe that also Proposition \ref{spreadoutpp} has
            a germ-grain counterpart.

            In this paper we have always fixed the distribution
            of the typical grain. However, it would be quite interesting
            to study the variability of empty space in germ-grain
            models for a fixed germ-process but variable grain distribution.
            For instance, one might compare models
            with equal expected volumes of the typical grains.
            To illustrate this task we give one simple example
            that is closely related to some of the results in \cite{SS80}.

\begin{example}\label{exBool}\rm Let $X_0$ be a random convex body
such that $E[V_d(X_0+K)]$ is finite for all compact $K\subset\R^d$.
Let $R$ and $\tilde R$ be positive random variables
with a finite $d$-th moment and assume that $X_0$ and
$R$ (resp.\ $X_0$ and $\tilde R$) are independent.  Consider two
Boolean model $Z$ and $\tilde Z$
based on the same Poisson process $N$ and
typical grains $R X_0$ and $\tilde R X_0$, respectively.
If
\begin{align*}
E[R^i]\ge E[\tilde R^i],\quad i=1,\ldots,d-1, 
\end{align*}
then (\ref{2.8}) holds for all structuring
elements $B$ such that $X_0$ and $B^*$ are a.s.\ in general relative
position, and all Borel sets $C\subset\R^d$. This follows from
\eqref{densBoolean} and the scaling property
$$
C_i(aX_0;B;\R^d\times C)=a^iC_i(X_0;B;\R^d\times C),\quad a> 0,
$$
see e.g.\ \cite{Schneider1993} for the Euclidean case $B=B(1)$.
\end{example}

\vspace{0.4cm} \noindent {\bf Acknowledgement}

\noindent
The first author is very grateful for the hospitality of the 
Mathematical Institute of the University of Wroc{\l}aw.
Large parts of this paper were written while he was visiting
the Institute.


\begin{thebibliography}{99}


%\bibitem[1994a]{BG94} A. Baddeley and R.D. Gill,
\bibitem{BG94} A. Baddeley and R.D. Gill (1994).
The empty space hazard of spatial pattern,
Preprint 845, University Utrecht, Department of Mathematics. 

%\bibitem[1997]{BG} A. Baddeley and R.D. Gill,
\bibitem{BG} A. Baddeley and R.D. Gill (1997).
\newblock Kaplan-Meier estimators of distance distributions for spatial point
processes, {\sl The Annals of Statistics} {\bf 25}, 263-292.


%\bibitem[2009]{blaszcz} B. B{\l}aszczyszyn and D. Yogeshwaran,
%\bibitem{blaszcz} B. B{\l}aszczyszyn and D. Yogeshwaran,
%Directionally convex ordering of random measures, shot noise
%fields, and some applications to wireless communications.
%{\sl Advances in Applied Probability} {\bf 41}, 623-646 (2009).

%\bibitem[2007]{BordTor07}
\bibitem{BordTor07}
C. Bordenave and G.L. Torrisi (2007).
Large deviations of Poisson cluster processes.
{\sl Stochastic Models} {\bf 23}, 593-625.

\bibitem{DaVJ03}
D. Daley and D. Vere-Jones (2003).
{\em An Introduction to the Theory of Point Processes}, 2nd edition.
Springer, New York.

\bibitem{Esary67}
J.D. Esary, F. Proschan and D. Walkaup (1967).
\newblock Association of random variables with applications.
{\sl Annals of Mathematical Statistics} {\bf 38}, 1466-1474.


\bibitem{Hall88}
P. Hall (1988).
{\em Introduction to the Theory of Coverage Processes},  Wiley, New York.

%\bibitem[1999]{BG99} Hansen, A. J. Baddeley, and R. D. Gill .
\bibitem{BG99} M. B. Hansen, A. J. Baddeley, and R. D. Gill (1999).
\newblock First contact distributions for spatial patterns:
regularity and estimation.
{\sl Advances in Applied Probability} {\bf 31}, 15-33.


\bibitem{H92} L. Heinrich (1992).
\newblock On existence and mixing properties of germ-grain models,
{\sl Statistics} {\bf 23}, 271-286.

\bibitem{Hl00}
D. Hug and G. Last (2000).
\newblock On support measures
in Minkowski spaces and contact distributions in stochastic
geometry. {\sl Annals of Probability}, {\bf 37}, 796-850.


\bibitem{HLW02a}
D. Hug, G. Last,  and W. Weil (2002).
\newblock Generalized contact distributions of inhomogeneous Boolean models.
{\sl Advances in Applied Probability} {\bf 34}, 21-47.

%\bibitem[2002b]{HLW02b}
\bibitem{HLW02b}
D. Hug, G. Last, and W. Weil (2002).
\newblock A survey on contact distributions.
Statistical Physics and Spatial Statistics, Lecture Notes in Physics {\bf 600},
317-357, {\sl Morphology of Condensed Matter, Physics and Geometry of
Spatially Complex Systems}, ed. by K. Mecke and D. Stoyan, Springer, Berlin.

%\bibitem[1999]{LH99} G. Last and M. Holtmann,
\bibitem{LH99} G. Last and M. Holtmann (1999).
\newblock On the empty space function of some germ-grain models, {\sl Pattern
Recognition}, {\bf 32}, 1587-1600.

\bibitem{LaScha89}
G. Last and  R. Schassberger (1998). On the distribution of
the  spherical contact vector of   stationary grain models.
{\sl Advances in Applied Probability} {\bf 30}, 36-52.


\bibitem{LiBa96}
M.N.M. van Lieshout and A. Baddeley (1996).  A nonparametric
measure of spatial interaction in point patterns,
{\sl Statistica  Neerlandica} {\bf 50}, 344-361.


\bibitem{Mo05}
I. Molchanov (2005).
\newblock{\em Theory of Random Sets},
\newblock Springer-Verlag, London.

%\bibitem[2002]{MulSto}
\bibitem{MulSto}
A. M\"{u}ller and D. Stoyan (2002).
\newblock {\em Comparison Methods for Stochastic Models and Risks},
\newblock Wiley, New York.

%\bibitem[1993]{Schneider1993}
\bibitem{Schneider1993}
R.~Schneider (1993).
\newblock {\em {C}onvex Bodies: the Brunn-Minkowski Theory},
\newblock Encyclopedia of Mathematics and its Applications 44, Cambridge
  University Press, Cambridge.


%\bibitem[2008]{SW08}
\bibitem{SW08}
R. Schneider and W. Weil (2008).
\newblock {\em Stochastic and Integral Geometry}.
Springer, Berlin.

\bibitem{St83} D. Stoyan (1983).
\newblock{\it Comparison Methods for Queues and Other Stochastic Models.} 
Wiley.


%\bibitem[1995a]{SKM95} D. Stoyan, W.S. Kendall and J. Mecke,
\bibitem{SKM95} D. Stoyan, W.S. Kendall and J. Mecke (1995).
\newblock {\it Stochastic Geometry and Its
Applications.} Second Edition, Wiley, Chichester.

\bibitem{SS80} H. Stoyan and D. Stoyan (1980). On some partial orderings of random closed sets.
{\sl Math. Operationsforsch. Stat. Optimization} {\bf 11}, 145-154.

\bibitem{Szekli95} R. Szekli (1995).
\newblock{\it Stochastic Ordering and Dependence in Apllied Probability.} 
Lecture Notes in Statistics 97, Springer Verlag, New York.

%\bibitem[1993a]{WW93} W. Weil and J.A. Wieacker,
\bibitem{WW93} W. Weil and J.A. Wieacker (1993).
\newblock Stochastic geometry. in:
{\sl Handbook of Convex Geometry}, North-Holland,
Amsterdam, 1391-1438.



\end{thebibliography}
\end{document}